\documentclass{amsart}

\usepackage[colorlinks=true, urlcolor=blue,bookmarks=true,bookmarksopen=true,citecolor=blue,hypertex]{hyperref}
\usepackage{graphicx}
\usepackage{enumerate}

\usepackage{amscd}
\usepackage{amssymb}
\usepackage{amsxtra}
\usepackage{amsmath}
\usepackage{enumerate}
\usepackage{mathrsfs}

\usepackage{color}

\usepackage{amsmath,amsfonts,amssymb}
\usepackage{verbatim}
\usepackage[hmargin = 4cm,vmargin = 2.5cm]{geometry}
\usepackage{epsfig}
\usepackage{amsthm}
\usepackage{bbold}

\usepackage[T1]{fontenc}

\theoremstyle{plain}

\newtheorem{thm}{Theorem}

\newtheorem{lem}[thm]{Lemma}

\theoremstyle{definition}

\theoremstyle{remark}

\newtheorem{rmk}[thm]{Remark}

\newtheorem{qs}[thm]{Question}

\theoremstyle{plain}

\newcommand{\Qed}{\hfill \qedsymbol \medskip}

\newcommand{\Id}{{{\mathchoice {\rm 1\mskip-4mu l} {\rm 1\mskip-4mu l}
      {\rm 1\mskip-4.5mu l} {\rm 1\mskip-5mu l}}}}
\def\NN{\mathbb{N}}
\def\ZZ{\mathbb{Z}}

\def\RR{\mathbb{R}}

\def\Id{\mathbb{1}}

\def\dd{\mathrm{d}}
\DeclareMathOperator{\Ham}{Ham}

\DeclareMathOperator{\Cal}{Cal}
\DeclareMathOperator{\flux}{flux}

\DeclareMathOperator{\Area}{Area}

\definecolor{lgray}{rgb}{0.5,0.5,0.5}

\begin{document}

\pagestyle{headings}

\bibliographystyle{alphanum}

\title{Curve obstruction for autonomous diffeomorphisms on surfaces}

\date{\today} 

\author{Michael Khanevsky}
\address{Michael Khanevsky, Faculty of Mathematics,
Technion - Israel Institute of Technology
Haifa, Israel}
\email{khanev@technion.ac.il}

\begin{abstract}
Consider a Hamiltonian diffeomorphism $g$ on a surface. We describe several scenarios where a curve $L$ and its image $g(L)$ 
provide a simple evidence that $g$ is not autonomous.
\end{abstract}

\maketitle

\section{Introduction and results}

An autonomous Hamiltonian flow on a surface admits an invariant lamination $\mathcal{B}$ by the flow lines. 
Therefore, the time-$t$ map $g$ of such flow preserves $\mathcal{B}$. It follows that given a curve $L$, both $L$ and $g(L)$
have the same pattern of intersection with $\mathcal{B}$. 
However, not any pair of curves $L, g(L)$ admits exisence of such a lamination $\mathcal{B}$. 
This is a potentially rich and unexplored source of combinatorial constraints that an autonomous map must satisfy. 
The present article is a little step in this direction, we describe an elementary obstruction in terms of combinatorics
of intersections of $L$ with $g(L)$ that prevents $g$ from having a loop of fixed points. 
At the same time, there are many scenarios where topological or dynamical arguments ensure existence of such loops for every $g$ generated by an autonomous flow.
In this case we obtain a tool which shows that if $g$ complies with the obstruction, it cannot be autonomous.

Fix a compact orientable symplectic surface $\Sigma$ (possibly with boundary) and a simple curve $L \subset \Sigma$. We ask that $L$ is either closed or 
connects boundary points of $\Sigma$, and does not touch $\partial \Sigma$ except at endpoints. We pick an orientation for $L$, it induces orientaion
also for $g(L)$.
This article restricts to flows $\phi$ on $\Sigma$ that are tangent to $\partial \Sigma$ (e.g. flows supported in the interior of $\Sigma$). 
Such $\phi$ admit a time-$t$ map for all $t \in \RR$ and preserve the connected components of $\partial \Sigma$ 
(as a set, but $\phi$ are allowed to rotate each boundary component).
Let $c$ be a periodic trajectory of $\phi$ (that is, $c$ is homeomorphic to a circle). All $x \in c$ are periodic with the same period $T_c$.
$\rho_x = 1/T_c$ is the \emph{rotation number} of $x$. $c$ is a \emph{loop of fixed points} of $\phi^1$ whenever $\rho_x \in \NN$ for some (and hence all) $x \in c$.

We say that a family $\mathcal{M} \subset \Ham (\Sigma)$ satisfies a \emph{fixed loop property detectable by $L$} (in short notation, $\mathcal{F}_L$ property)
if every autonomous $g \in \mathcal{M}$ has a loop of fixed points which intersects $L$.
In Section~\ref{S:ex} we provide few scenarios where such property can be established by topological or dynamical considerations. 

\begin{figure}[!htbp] 
\begin{center}
\includegraphics[width=0.6\textwidth]{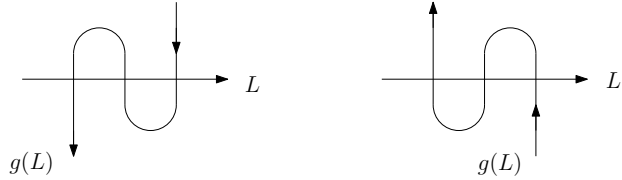}
\caption{Local obstruction}
\label{F:loc-obs}
\end{center}
\end{figure}

The main result is the following: 
\begin{thm} \label{T:main}
  Suppose $\mathcal{M} \subset \Ham (\Sigma)$ has the $\mathcal{F}_L$ property, let $g \in \mathcal{M}$.
  
  Assume $L \pitchfork g(L)$ away from a neighborhood of $\partial \Sigma$ and there is a system of disjoint neighborhoods $\{U_i\}$ for transverse intersections in $L \cap g(L)$, 
  such that in each neighborhood $(L \cup g (L)) \cap U_i$ is diffeomorphic to one of the configurations on Figure~\ref{F:loc-obs}. Then $g$ is non-autonomous.
\end{thm}

 This result can be compared to ~\cite{Kh:non-auton-curves} which establishes a different kind of obstructions for a Hamiltonian $g$ to be autonomous. 
 Similar to this article, the obstructions there can be described in terms of $L$ and its image $g(L)$. However, ~\cite{Kh:non-auton-curves} relies on more involved tools 
 (various quasimorphisms) which reflect global dynamics of $g$, while here we use an elementary argument which is purely local near the intersection points.

\section{The obstruction}

Let $\phi$ be an autonomous Hamiltonian flow on $\Sigma$.
We call a closed trajectory $c$ \emph{regular} with respect to $L$ if $c \pitchfork L$ and $x$ is not a critical point for $\rho$.
This implies that near $c$ one may find periodic trajectories $c^{-}$ and $c^+$ with periods $T_{c^-} < T_c < T_{c^+}$.

Given a loop of fixed points $c$ of $\phi^1$ such that $c \cap L \neq \emptyset$, one can perturb $\phi^1$ in $C^1$ topology 
into an autonomous Hamiltonian which has a regular loop of fixed points near $c$: if $c \cap L$ is not transverse, conjugate $\phi$ with a $C^\infty$-small diffeomorphism 
$h$ such that $h(c) \cap L \neq \emptyset$ and $h (c) \pitchfork L$. $h(c)$ is a loop of fixed points for the perturbed flow and has the same rotation number as $c$ with respect to $\phi$. 
If $T_c$ is a critical value, one can accelerate (or slow down) the flow inside nearby trajectories on either side of $c$ and destroy this condition.

\medskip

Suppose that $c$ is regular for $g = \phi^1$. Note, that since $c$ is fixed pointwise by $g$, $c \cap L = c \cap g (L)$. That is, $c$ 
is allowed to intersect $L \cup g(L)$ only at $L \cap g(L)$. Moreover, as $g$ is orientation-preserving, $L$ and $g (L)$ induce the same coorientation at their transverse 
intersections with $c$. In simple words, at an intersection point $x \in c \cap L \cap g(L)$, both $L$ and $g(L)$ point to the same side of $c$. 
If we pick a coordinate chart around $p \in L \cap g (L)$
which sends $L$ to the $x$-axis and $g(L)$ to the $y$-axis (respecting orientations of curves), then $c$ may cross the intersection either from the second quadrant to the 
fourth or vice versa. It cannot connect the first and the third quadrants as in this case $L$ and $g(L)$ have opposite coorientations.

\begin{figure}[!htbp] 
\begin{center}
\includegraphics[width=0.6\textwidth]{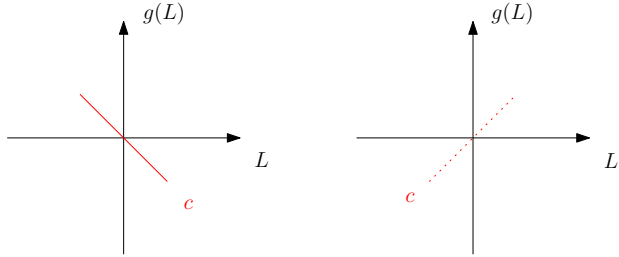}
\caption{$L, g(L)$ have the same coorientation with respect to $c$ (left) and opposite coorientation (right)}
\label{F:quad}
\end{center} 
\end{figure}

We describe below a local obstruction to have a \emph{regular} loop of fixed points $c$ passing through $x \in L \cap g(L)$. 
We provide two simple proofs for this phenomenon, one purely combinatorial and another having dynamical nature.

\begin{figure}[!htbp] 
\begin{center}
\includegraphics[width=0.8\textwidth]{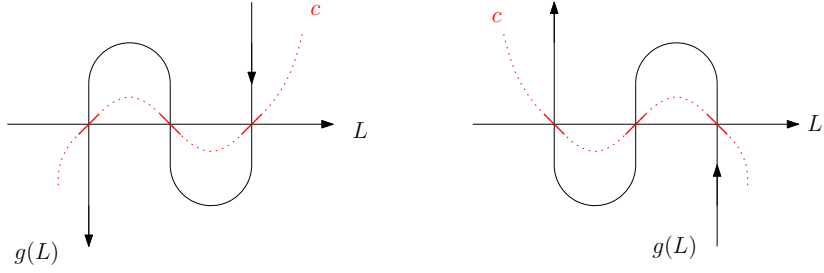}
\caption{$c$ crossing $L$}
\label{F:loc-cross}
\end{center} 
\end{figure}

\begin{lem}
Suppose there is a coordinate chart $U$ which contains three intersection points of $L \cap g(L)$ connected as it is shown in Figure ~\ref{F:loc-obs}.
Then none of these intersection points lies on a regular loop of fixed points of $g$.
\end{lem}

Proof 1:
 Assume by contradiction that $c$ is such a loop.
 In order to keep coorientation of $L$ consistent with that of $g(L)$, $c$ has to intersect the curves in direction indicated by 
 solid red lines in Figure ~\ref{F:loc-cross}. Moreover, if $c$ crosses from one side of $L$ to another and avoids $L \cup g (L)$ except at their intersection points, 
 it has to connect all three intersection points as described by the dotted line. 

  We look at the order of these three intersection points: the order along $L$ does not match that along $g(L)$ 
  (even if $L$ is a closed curve and we consider the circular order, it is still different). But if the points are fixed for $g$, their order must be preserved.
  Indeed, let $\gamma$ be the oriented arc of $L$ which connects two intersection points $x, y \in L \cap g(L)$. Then $g(\gamma)$ is the oriented arc
  in $g(L)$ which connects $x$ to $y$. A chain of three intersection points $x, y, z$ connected by disjoint arcs $\gamma_1, \gamma_2$ is sent by $g$ to the chain
  $x, y, z$ connected by $g(\gamma_1), g(\gamma_2)$ -- the order of $x, y, z$ along $g(L)$ is the same as along $L$.
\Qed

Proof 2:
 Assume by contradiction that $c$ is such a loop.
 As in Proof 1, coorientations constrain $c$ to intersect the curves as indicated in the figure. 
 
 Given a fixed point $x \in c$ and a short oriented curve $\gamma$ through $x$ which is transverse to $c$, we call $x$ \emph{positive} 
 if the angle between $\gamma$ and $g (\gamma)$ at $x$ is in counterclockwise direction. When $c$ is a regular loop of fixed points, this angle cannot be zero.
 Note that positivity does not depend on the choice of $\gamma$: pick $v_c, v_\gamma \in T_x\Sigma$ - non-zero vectors tangent to $c$, $\gamma$, respectively,
 such that $(v_c, v_\gamma)$ is a basis with negative orientation. The linearization $Dg_x :T_x \Sigma \to T_x\Sigma$ is represented in this basis
 by the matrix $\left( \begin{array}{cc} 1 & \alpha \\ 0 & 1 \end{array}\right)$. Indeed, $c$ is fixed pointwise, hence $v_c$ is an eigenvector with the eigenvalue $1$. 
 $g$ is Hamiltonian implies it is area-preserving hence the second eigenvalue must be $1$ as well. 
 $c$ is regular means $Dg_x$ cannot be the identity map, so $\alpha \neq 0$.
 The sign of $x$ coincides with the sign of $\alpha$ and $\alpha$ depends continuously on $\gamma'(x)$.
 
 Note as well that the sign is constant at all points $x \in c$. However, if $c$ had intersected $L \cap g(L)$ as indicated, it would have had
 alternating signs at the three intersection points - a contradiction.
\Qed

\begin{rmk}
  Our obstruction is stable under $C^1$-small perturbations of $g$ (such perturbations do not change the combinatorics of $L \cap g(L)$). 
  It follows that the obstruction prevents existence of general loops of fixed points for $g$ (not only regular ones) that intersect $L$ at $x$: 
  assume by contradiction that $g$ is autonomous with a non-regular loop of fixed points $c$. 
  We perturb $g$ into $g'$ as it was described in the beginning of this section. $g'$ has a regular loop of fixed points $c'$ -- but the obstruction persists also for $g'$. 
  A contradiction.
\end{rmk}

If all intersection points $L \cap g(L)$ satisfy the obstruction, $g$ does not have a loop of fixed points which intersects $L$. 
This immediately implies Theorem~\ref{T:main}.

In the next section we describe several examples of $(\Sigma, L)$ and families $\mathcal{M} \subset \Ham(\Sigma)$ where one may establish 
the $\mathcal{F}_L$ property. In these scenarios our obstruction provides both an easy criterion to show 
non-autonomy of $g$ and a method to perturb $g (L)$ into $L'$ so that any $\{ h \in \mathcal{M} \, | \, h(L) = L' \}$
is not autonomous. Namely, go over all transverse intersection points $L \cap g(L)$ and replace a small arc of $g(L)$ around the intersection by a ``snake'' described in Figure~\ref{F:loc-obs}
(pattern to the left if $g (L)$ crosses $L$ from the left to the right or pattern to the right for the opposite direction). 
This perturbation may be realized by a Hamiltonian $h$ supported near the intersection point such that $h$ has arbitrarily small $C^0$ and Hofer norms.

\begin{figure}[!htbp] 
\begin{center}
\includegraphics[width=0.7\textwidth]{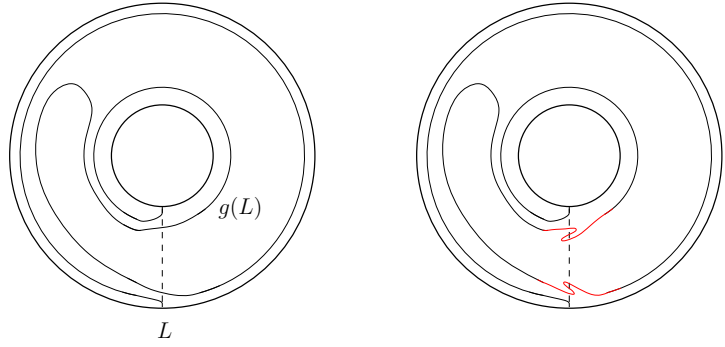}
\caption{Unobstructed (left) and and obstructed (right) intersections after a local perturbation}
\label{F:obs-pert}
\end{center} 
\end{figure}

\begin{qs}
  A standard method to verify that $g$ is not autonomous is to show that it has positive entropy. That can often be done locally
  or semi-locally as it is enough to show existence of horseshoe dynamics (e.g. via transverse homoclinic points). 
  However, the obstruction described above does not use the full data of $g$: for example, it is difficult to say something certain 
  about the fixed points looking just at $L$ and $g(L)$. Is there a common background for these methods? 
  Does the obstruction imply positive entropy for $g$ (rather than just showing that $g$ is not autonomous)?
\end{qs}

\section{Examples} \label{S:ex}

\subsection{The annulus}

Let $A = S^1 \times [0, 1]$ (we set $S^1 = \RR / \ZZ$) equipped with the area form $\omega = \dd \theta \wedge \dd s$ so that $A$ has total area $1$.
For an autonomous flow $\phi$ on $A$ we may consider rotation numbers $r^\phi_A: A \to \RR$ with respect to the $S^1$ coordinate of $A$. Looking at the Reeb graph
of the Hamiltonian function which generates $\phi$, it is easy to show that $r_A (x)$ exists for all $x \in A$ and is continuous.
Moreover, for periodic points $x$ with noncontractible orbits $r_A (x) = \pm \rho_x$ which was defined earlier (if the orbit of $x$ is either contractible
or not closed, $r_A (x) = 0$).

In all examples below we set $L = \{ 0 \} \times [0, 1]$. We describe families $\mathcal{M} \subset \Ham (A)$ which guarantee existence of $x \in A$ with 
$r_A (x) \in \ZZ \setminus \{0\}$ for all autonomous $g \in \mathcal{M}$. That implies that the orbit $c$ of $x$ is a closed non-contractible loop with 
$\rho_x \in \ZZ$. As all non-contractible loops in $A$ intersect $L$, these $\mathcal{M}$ satisfy the fixed loop property detectable by $L$.
(We have to exclude points with $r_A (x) = 0$ from our scope of interest as they may be non-periodic or belong to orbits which avoid $L$.) 

\subsubsection{Example}
given $g \in \Ham (A)$ let $\tilde{g} \in \Ham (\RR \times [0, 1])$ be its lift to the universal cover of $A$.
\[
  \mathcal{M} = \{g \in \Ham (A) \, | \, \pi_\RR (\tilde{g} (0, 1)) - \pi_\RR (\tilde{g} (0, 0)) > 4\}.
\]
Note that while the projection $\pi_\RR (\tilde{g} (0, 1))$ depends on the lift, the difference 
$\pi_\RR (\tilde{g} (0, 1)) - \pi_\RR (\tilde{g} (0, 0))$ is determined by $g$.
Moreover, it may be described in the following way: let $\widetilde{g(L)}$ be a lift of $g(L)$ to $\RR \times [0, 1]$.
$$\pi_\RR (\tilde{g} (0, 1)) - \pi_\RR (\tilde{g} (0, 0)) = \pi_\RR (\tilde{y})) - \pi_\RR (\tilde{x})$$ where $\tilde{x}, \tilde{y}$ are the endpoints of 
$\widetilde{g(L)}$. That is, the family $\mathcal{M}$ is determined by the curve $g(L)$.

\medskip

It is well-known that for an autonomous flow $\varphi : S^1 \to S^1$ the rotation number $\rho^\varphi$ satisfies 
$\widetilde{\varphi}^1 (x) - x - 1 < \rho^\varphi < \widetilde{\varphi}^1 (x) - x + 1$ 
where $\widetilde{\varphi} : \RR \to \RR$ is the lift of $\varphi$ to the universal cover such that $\widetilde{\varphi}^0 = \Id$ and $x \in \RR$ is arbitrary.
Suppose $g \in \mathcal{M}$ is the time-$1$ map of an autonomous Hamiltonian flow $\phi$. $\phi$ restricts to a flow on the boundary components so
\[
  r^\phi_A (0, 1) - r^\phi_A (0, 0) >  (\pi_\RR (\tilde{g} (0, 1)) - 1) - (\pi_\RR (\tilde{g} (0, 0)) + 1) > 2.
\]
By continuity $r^\phi_A$ has at least two different integer values in its image, so at least one value in $\ZZ \setminus \{0\}$. It follows that
$\mathcal{M}$ satisfies the $\mathcal{F}_L$ property.

\begin{rmk}
  The rotation number $r^\phi_A$ may depend on the isotopy $\phi$ rather than on its time-$1$ map $g$, different choices of isotopy may give rise to 
  rotation numbers shifted by an integer. We have shown above that there exists an orbit with non-zero integer rotation number for every choice of 
  autonomous isotopy that generates $g$.
\end{rmk}

\subsubsection{Example} \label{SSS:a_flux}
In this example we restrict to $\Ham' (A) \subset \Ham (A)$ generated by Hamiltonian functions $H_t : A \to \RR$ which are locally constant near 
$\partial A$ (that is, the Hamiltonian flow is stationary near the boundary, but the flux through $L$ might be not zero).
\[
  \mathcal{M} = \{g \in \Ham' (A) \, | \, \flux_L (g) > 1\}.
\]
The flux may be computed by $$\flux_L (g) = \int H_t (1, 0) \dd t - \int H_t (0, 0) \dd t$$ where $H_t$ is a Hamitonian function which generates $g$.
Alternatively, $\flux_L (g)$ is equal to the area bounded between $\widetilde{L}$ and $\tilde{g} (\widetilde{L})$ where $\widetilde{L}$ is a lift of $L$ to the universal cover of $A$
and $\tilde{g}$ is the lift of $g$ which restricts to the identity near the boundary.

In this setting $r_A (x)$ is well defined for all autonomous $g$ (does not depend on the isotopy $\phi \subset \Ham' (A)$). Moreover, for autonomous $g$
\[
  \flux_L (g) = \int_A r_A (x) \omega > 1
\]
(this identity is easy to verify using action-angle coordinates near periodic orbits of $\phi$).
Therefore there exists $x \in A$ with $r_A (x) \geq 1$. At the same time $r_A = 0$ near $\partial A$. Continuity implies that $r_A$ 
has value $1$ in its image and the $\mathcal{F}_L$ property follows.

\subsubsection{Example}  \label{SSS:a_calabi}
We restrict to $\Ham_c (A) \subset \Ham (A)$ generated by Hamiltonian functions $H_t : A \to \RR$ which are compactly supported in the interior of $A$.
(that is, the Hamiltonian flow is stationary near the boundary, and the flux through $L$ is zero).

\medskip

Let $F_t : \Sigma \to \RR$, $t \in [0, 1]$ be a time-dependent smooth function with compact support in the interior of $\Sigma$. We define
$\widetilde{\Cal} (F_t) = \int_0^1 \left( \int_\Sigma F_t \omega \right) \dd t$. If the symplectic form $\omega$ is exact 
(this is the case for an annulus or a disk),
$\widetilde{\Cal}$ descends to a homomorphism $\Cal_\Sigma: \Ham(\Sigma) \to \RR$ which is called 
the \emph{Calabi homomorphism}.

\medskip

We also need the \emph{Calabi quasimorphism} by Entov-Polterovich (see ~\cite{En-Po:calqm}). The authors construct a homogeneous 
quasimorphism $\Cal_{S^2} : \Ham (S^2) \to \RR$. It has many wonderful properties related to the group structure
and geometry of $\Ham (S^2)$, here we will state just the following one: for an autonomous Hamiltonian $f : S^2 \to S^2$,
\[
	\Cal_{S^2}(f) = \int_{S^2} F \omega - \Area (S^2) \cdot F(X),
\]
where $F$ is a generating function for $f$ and $X$ is a certain level set of $F$ (the \emph{median}, consider ~\cite{En-Po:calqm} for details).

We embed $A$ into a sphere $S^2_{a,b}$ of area $1 + a + b$ by gluing a disk of area $a$ to $S^1 \times \{0\}$ and a disk of area $b$ to $S^1 \times \{1\}$.
Denote this embedding by $i_{a, b} : A \to S^2_{a,b}$.
A Hamiltonian function $F:A \to \RR$ supported in the interior of $A$ extends by zeros to a function $S^2 \to \RR$. This allows us to pull $\Cal_{S^2}$
back to a quasimorphism $i_{a,b}^* \Cal_{S^2_{a,b}} : \Ham_c (A) \to \RR$.
Let $$r_{a,b} = \frac{1}{1+a+b} \cdot \left( \Cal_A - i_{a,b}^* \Cal_{S^2_{a,b}} \right)$$ be the normalized difference between the 
Calabi homomorphism on $A$ and the pullback of the Calabi quasimorphism. Using methods from ~\cite{En-Po-Py:qm-continuity} one shows that
$r_{a,b}$ vanishes on Hamiltonians supported in a small disk hence is $C^0$-continuous.

Fix a parameter $h \in [0, 1]$. Let $G: A \to \RR$ be a Hamiltonian function, $g$ its time-$1$ map and set $a = 1, b = 1 + 2h$. 
For a generic $G$, $r_{1, 1+2h} (g) = G(X_h)$ where $X_h$ is the level set of $G$ which is sent to a median of $i_{1, 1+2h, *} (G)$ on $S^2$.
One may show that $X_h$ is not contractible in $A$ and $X_h$ together with $S^1 \times \{0\}$ bound a subannulus $A_g^h$ of area at most $h$.
The flux of $g$ in $A_g^h$ is precisely $G(X_h) = r_{1, 1+2h} (g)$. From here we continue as in the previous example.

\medskip

\[
  \mathcal{M}_h = \{g \in \Ham_c (A) \, | \, r_{1, 1+2h} (g) > h\}.
\]

As before, for generic autonomous $g$
\[
  \flux_{A_g^h} (g) = \int_{A_g^h} r_A (x) \omega > h
\]
hence there exists $x \in A$ with $r_A (x) = 1$. $r_{1, 1+2h}$ is $C^0$-continuous, so applying a perturbation we may drop the requirement that 
$g$ is generic. The $\mathcal{F}_L$ property follows.

\begin{rmk}
  Using tools described in ~\cite{Kh:non-auton-curves}, one may show that the condition $r_{1, 1+2h} (g) > h$ can be 
  determined (up to a bounded defect) by the pair $L, g(L)$.
\end{rmk}

\subsection{The disk} \label{SS:disk}
Let $D$ be the unit disk in $\RR^2$ equipped with the area form $\omega = \frac{1}{\pi} \dd x \wedge \dd y$ so that $D$ has total area $1$.
$L$ is the diameter $[-1, 1] \times \{0\}$.
We use a similar construction to Example ~\ref{SSS:a_calabi}

Embed $D$ into a sphere $S^2_{a}$ of area $1 + a$ by gluing a disk of area $a$ to $\partial D$.
Denote this embedding by $i_{a} : D \to S^2_{a}$.
Let $$r_{a} = \frac{1}{1+a} \cdot \left( \Cal_D - i_{a}^* \Cal_{S^2_{a}} \right)$$ be the normalized difference between the 
Calabi homomorphism on $D$ and the pullback of the Calabi quasimorphism of $S^2_{a}$.

Let $G: D \to \RR$ be a Hamiltonian function compactly supported in the interior of $D$, $g$ its time-$1$ map and set $a = 1 - 2h$ 
for a fixed parameter $h \in (0, 0.5)$. 
For a generic $G$, $r_{1-2h} (g) = G(X_h)$ where $X_h$ is the level set of $G$ which is sent to the median of $i_{1-2h} (G)$ on $S^2$.
One shows that if $G(X_h) \neq 0$, $X_h$ together with $\partial D$ bound a subannulus $A_g^h$ of area at most $h$.
The flux of $g$ in $A_g^h$ is precisely $G(X_h) = r_{1-2h} (g)$. From here we continue as before.

\medskip

\[
  \mathcal{M}_h = \{g \in \Ham_c (D) \, | \, r_{1-2h} (g) > h\}.
\]

Once again, for generic autonomous $g$
\[
  \flux_{A_g^h} (g) = \int_{A_g^h} r_{A_g^h} (x) \omega > h
\]
hence there exists $x \in A$ with $r_{A_g^h} (x) = 1$. The orbit $c$ of $x$ is a non-contractible loop in $A_g^h$ hence $c$ encircles
$D \setminus A_g^h$. Therefore the region bounded by $c$ is a tolopogical disk of area at least $1-h > 1/2$. By area considerations, $c$ must intersect $L$.
$r_{1-2h}$ is $C^0$-continuous, so applying a perturbation we may drop the requirement that 
$g$ is generic. 
The $\mathcal{F}_L$ property follows.

\begin{rmk}
  Using the argument from ~\cite{Kh:diam}, one shows that the condition $r_{1-2h} (g) > h$ can be determined 
  (up to a bounded defect) by the pair $L, g(L)$.
\end{rmk}

%
%
%

\subsection{Final remarks}
\begin{rmk}
  In close analogues of the examples above the family $\mathcal{M}$ can be described in terms of the curve $g(L)$ rather than the map $g$ itself. That is, we have constructed 
  an obstruction that depends only on the pair of curves $L$ and $g(L)$.
\end{rmk}

\begin{rmk}
  In all examples except for the first one, we may refine the obstruction. Instead of blocking all possible intersections of loops of fixed points $c$
  with $L$ it is enough to obstruct only those intersection points $x \in L \cap g(L)$ that have Maslov index $2$, as these are the 
  intersection points that arise when we have rotation number $r_A(x) = 1$. 
\end{rmk}

\bibliography{bibliography}

\end{document}